\documentclass[11pt]{article}
  \usepackage{latexsym}
  \usepackage{amssymb}

\usepackage{a4wide}
   \usepackage{geometry}
   \usepackage{lscape}
   \usepackage{rotating}
  \makeatletter
  \def\section{\@startsection {section}{1}{\z@}{-3.5ex plus -1ex minus
   -.2ex}{2.3ex plus .2ex}{\large\bf}}
  \def\subsection{\@startsection{subsection}{2}{\z@}{-3.25ex plus -1ex minus
   -.2ex}{1.5ex plus .2ex}{\normalsize\bf}}
  \makeatother

  \makeatletter
  \@addtoreset{equation}{section}

 \newenvironment{prof}[1][Proof]{\textbf{#1.} }{\ \rule{0.5em}{0.5em}}

\newtheorem{theorem}{Theorem}[section]
\newtheorem{definition}{Definition}[section]

\newtheorem{example}{Example}[section]

  \makeatother

\begin{document}
%

\begin{center}

\vspace{1cm}

\vskip .5in
 {\Large\bf General second order  conditions for extrema of functionals}

\vspace{10pt}

Mahouton Norbert Hounkonnou$^{1,\dag}$ and Pascal Dkengne Sielenou$^{1,*}$

$^{1}${\em University of Abomey-Calavi,\\ International Chair in Mathematical Physics
and Applications}\\
{\em (ICMPA--UNESCO Chair), 072 B.P. 50  Cotonou, Republic of Benin}\\

E-mails:   $^{\dag}$norbert.hounkonnou@cimpa.uac.bj, \quad $^{*}$sielenou$\_$alain@yahoo.fr
 \vspace{0.5 cm}
\begin{abstract}
We  prove both necessary and sufficient second order conditions
 of extrema for  variational problems involving any  higher order continuously  twice differentiable
 Lagrangians with multi-valued dependent functions of several variables.
 Our analysis is performed in the framework of the finite dimensional total jet space.
\end{abstract}
\end{center}
\makeatother

\section{Introduction}

The calculus of variations encompasses a very broad range of mathematical applications.
The methods of variational analysis can be applied to an enormous variety of physical systems,
whose equilibrium configurations  minimize or maximize a suitable functional,
which typically represents the potential energy of the system.    The
critical functions are characterized as solutions to a system of partial differential equations,
known as the Euler-Lagrange equations derived from the variational principle.
Each solution to the boundary value problem specified by the Euler-Lagrange
equations subject to appropriate boundary conditions is thus a candidate extremizer for the variational
problem.  In many applications, the Euler-Lagrange boundary value problem suffices to
single out the physically relevant solutions, and one needs not press onto the considerably
more difficult second variation.

In general, the solutions to the Euler-Lagrange boundary value problem are critical
functions for the variational problem, and hence include all (smooth) local and global extremizers.
The determination of which solutions are genuine minima or maxima requires a further analysis
of the positivity properties of the second variation.
Indeed, as stated  in \cite{r09}, a complete analysis of the positive definiteness of the second
variation of multi-dimensional variational problems is quite complicated, and still awaits
a completely satisfactory resolution!
This  is thus a reason for which  second order conditions of extrema is  customary  established only for at most two-dimensional variational
problems involving first order Lagrangians \cite{r01,r02,r03}.
The aim  of this paper is to formulate the second order extrema conditions for multi-dimensional
 variational problems with any higher order Lagrangians.

\section{Preliminaries: basic definitions, theorems and notations}
This section, mainly based on \cite{r10, r13}, addresses relevant  definitions, theorems and notations playing
a central role in the  calculus of variations.

Let $x=(x^1,\ldots,x^n)$ and $u=(u^1,\ldots,u^m),$ $u=u(x),$
 and $X\times U^{(s)},$ the space whose coordinates are
denoted by $(x,u^{(s)}),$ encompassing    the independent variables $x$,  the
dependent variables $u$ and their derivatives up to order $s$, $u^{(s)}$.
\begin{definition}\textbf{\emph{(Differential function)}}\\
A function $f$ defined on  $X\times U^{(s)}$ is called
  $s$-order differential function   if it is locally analytic, i.e., locally expandable in a Taylor series with respect to all arguments.
\end{definition}
\begin{definition}\emph{\textbf{(Total derivative operator)}}\\
 Let $f$  defined on  $X\times U^{(s)}$ be an $s$-order differential function.
 The total derivative of $f$ with respect to $x^i$ is defined by:
 $$
   D_{x^i}f=\frac{\partial f}{\partial
   x^i}+\sum_{j=1}^m\sum_{k_1+\cdots+k_n=0}^{s}u^j_{k_1x^1\cdots (k_i+1)x^i\cdots k_nx^n}\frac{\partial f}{\partial u^j_{k_1x^1\cdots
   k_nx^n}},
 $$
where
$$
 u^j_{k_1x^1\cdots
   k_nx^n}=\frac{\partial^{k_1+\cdots+k_n}u^j}{(\partial x^1)^{k_1}\cdots (\partial
   x^n)^{k_n}},
$$
\end{definition}
\begin{definition}\textbf{\emph{(Zeroth-Euler operator)}}\\
Let $f$ defined on $X\times U^{(s)}$ be an $s$-order  differential function.
 The zeroth-Euler operator (also called  the variational derivative)
  of $f$ is given by
  \begin{equation}
  \frac{\delta}{\delta u}f=\left(\frac{\delta}{\delta u^1}f,\cdots,\frac{\delta}{\delta u^m}f\right),
  \end{equation}
where for $j=1,\cdots,m$
\begin{equation}
\frac{\delta}{\delta u^j}f=\sum_{k_1+ \cdots  +k_n=0}^{s}(-D_{x^1})^{k_1}\cdots(-D_{x^n})^{k_n}\frac{\partial
f}{\partial u^j_{k_1x^1\cdots k_nx^n}}.
\end{equation}
\end{definition}
A variational problem consists in finding extrema of a functional $J$ defined by
\begin{equation}\label{eq6}
J[u]=\int_{\Omega}L\left(x,u^{(s)}(x)\right)dx,
\end{equation}
where $\Omega$ is a connected open subset of $X$ and $L$ defined on $X\times U^{(s)}$ is an  $s$-order differential function
 called the Lagrangian of the variational problem $J.$
 In general, a functional is a mapping that assigns to each element in some
function space a real number, and a variational problem amounts to searching for functions which are an extremum (minimum, maximum) or saddle points of a given functional.
\begin{theorem}
Let $u$ be an extremal of $J,$ then $u$ satisfies the Euler-Lagrange equations
\begin{equation}\label{eq7}
 \frac{\delta}{\delta u^j} L\left(x,u^{(s)}(x)\right)  =0,\quad j=1,\cdots,m.
\end{equation}
\end{theorem}


\section{Necessary and sufficient conditions for extrema of functionals}

In this section, we propose a definition of the  total jet space and
show that in this framework,  the shape of second order
conditions for extrema of functional whose Lagrangian includes multi-valued dependent functions of several variables remains the same 
as that of the  second variation for a functional   with first order one-dimensional scalar valued Lagrangian.

Consider $X$, an $n$-dimensional independent variables space, and $U=\bigotimes_{j=1}^{m}U^j,$
 an $m$-dimensional dependent variables space.
 Let $x=(x^1,\cdots,x^n)\in X$ and $u=(u^1,\cdots,u^m)\in U$ with $u^j\in U^j.$
  We define  the jet-space $U^{(s)}$ as:
  \begin{equation}
   U^{(s)}:=\bigotimes_{j=1}^{m}\left(\bigotimes_{l=0}^{s}U^j_{(l)}  \right),
  \end{equation}
 where $U^j_{(l)}$ is the set of all $p_l\equiv \left(\begin{array}{c}{n+l-1}\\{l}  \end{array}\right)$ distinct $l$-th order
 partial derivatives of $u^j.$ We denote by
  $u^j_{(k)}$ the $p_k$-tuple of all $k$-order partial derivatives of $u^j$.
 The   $u^j_{(k)}$ vector components are recursively obtained  as follows:
 \begin{itemize}
 \item[i)] $u^j_{(0)}=u^j$ and $u^j_{(1)}=(u^j_{x^1},u^j_{x^2},\cdots,u^j_{x^n}).$
 \item[ii)] Assume that  $u^j_{(k)}$ is known and form the tuples
 $$\widetilde{u}^j_{(k+1)}(l)=\left(\frac{\partial}{\partial x^1}u^j_{(k)}[l], \frac{\partial}{\partial x^2}u^j_{(k)}[l],\cdots,\frac{\partial}{\partial x^n}u^j_{(k)}[l]   \right),\quad l=1,2,\cdots,p_k;   $$
 $$\widetilde{u}^j_{(k+1)}=\left( \widetilde{u}^j_{(k+1)}(1),\widetilde{u}^j_{(k+1)}(2),\cdots,\widetilde{u}^j_{(k+1)}(p_k)   \right),   $$
 where $u^j_{(k)}[l]$ is the $l$-th component of the vector $u^j_{(k)}.$
 \item[iii)] The tuple $u^j_{(k+1)}$ is obtained from the tuple $\widetilde{u}^j_{(k+1)}$ in such a way to exclude all further components already
written following the  vector components order.
 \end{itemize}
The simplest instances of the above situation are obtained when $n=2$ and $3$ as follows.
\begin{example}:
\begin{itemize}
 \item For
   $n=2$,  $x=(x^1,x^2)$ and we have:
  \begin{eqnarray}
   u^j_{(1)}=\left(u^j_{x^1},u^j_{x^2}\right),\nonumber 
  \end{eqnarray}
and
\begin{eqnarray}
 \widetilde{u}^j_{(2)}(1)&=&\left(\frac{\partial}{\partial x^1}u^j_{(1)}[1], \frac{\partial}{\partial x^2}u^j_{(1)}[1]   \right)
 =\left(u^j_{2x^1},u^j_{x^1x^2}\right),\nonumber\cr
 \widetilde{u}^j_{(2)}(2)&=&\left(\frac{\partial}{\partial x^1}u^j_{(1)}[2], \frac{\partial}{\partial x^2}u^j_{(1)}[2]   \right)
 =\left(u^j_{x^2x^1},u^j_{2x^2}\right),\nonumber\cr
\widetilde{u}^j_{(2)}&=&\left(\widetilde{u}^j_{(2)}(1),\widetilde{u}^j_{(2)}(2)  \right) 
=\left( u^j_{2x^1},u^j_{x^1x^2},  u^j_{x^2x^1},u^j_{2x^2} \right),\nonumber\cr
u^j_{(2)}&=&\left(u^j_{2x^1}, u^j_{x^1x^2},  \check{u}^j_{x^2x^1} ,u^j_{2x^2} \right)
 =\left(u^j_{2x^1},u^j_{x^1x^2},u^j_{2x^2}\right).\nonumber
\end{eqnarray}
  \item  For $n=3$, $x=(x^1,x^2,x^3)$ and we get
 \begin{eqnarray} u^j_{(2)}&=&\left(u^j_{2x^1},u^j_{x^1x^2},u^j_{x^1x^3},u^j_{2x^2},u^j_{x^2x^3},u^j_{2x^3}  \right),\nonumber\cr
u^j_{(3)}&=&\left(u^j_{3x^1},u^j_{2x^1x^2},u^j_{2x^1x^3},u^j_{x^12x^2},u^j_{x^1x^2x^3},u^j_{x^12x^3},u^j_{3x^2},u^j_{2x^2x^3},u^j_{x^22x^3},u^j_{3x^3}  \right),\nonumber\end{eqnarray}
for $k=2$ and $k=3$, respectively.
\end{itemize}
\end{example}
  An element   $u^{(s)}$ in the jet-space $U^{(s)}$ is the $m(1+p_1+p_2+\cdots+p_s)=m \left(\begin{array}{c}{n+s}\\{s}  \end{array}\right)$-tuple defined by
$$u^{(s)}=\left(u^1_{(0)},u^1_{(1)},\cdots,u^1_{(s)}, u^2_{(0)},u^2_{(1)},\cdots,u^2_{(s)},\cdots,u^m_{(0)},u^m_{(1)},\cdots,u^m_{(s)}  \right).$$
Naturally, we have
$$ u^{j(s)}=\left(u^j_{(0)}, u^j_{(1)},u^j_{(2)}, \cdots , u^j_{(s)}\right),  \quad j=1,2, \cdots  ,m.$$
\begin{definition}
Consider an $s$-order variational problem defined by the functional
\begin{equation}\label{eqn0}
J[u]=\int_{\Omega}L\left(x,u^{(s)}(x)\right)dx,
\end{equation}
where  $\Omega$ is an open subset of $\mathbb{R}^n.$ Then, we
 define an $m\times m$-matrix $A$ of second order partial derivatives of $L$ by:
\begin{equation}\label{eqac}
A=\left[ A^{jj'} \right]_{1\leq j,\,j'\leq m}\quad \emph{\emph{with}}\quad A^{jj'}=\left[ A^{jj'}_{kk'} \right]_{0\leq k,\,k'\leq s},
\end{equation}
 \begin{equation}
A^{jj'}_{kk'}=\left[ \frac{\partial^2L}{\partial u^j_{(k)}[h]\partial u^{j'}_{(k')}[h']} \right]_{_{1\leq h\leq p_k,\,}^{1\leq h'\leq p_{k'}}}.
\end{equation}
\end{definition}
\begin{example}
 Let us construct the martrix $A$ for particular values of the integers $m,\,n,\,s.$  If $m=s=2$, then
 $$A=\left[\begin{array}{cc} A^{11}& A^{12}\\A^{21}&A^{22}\end{array}\right] \quad \emph{\emph{with}}   \quad A^{jj'}=\left[\begin{array}{ccc} A^{jj'}_{00}& A^{jj'}_{01}&A^{jj'}_{02}\\
 A^{jj'}_{10}&A^{jj'}_{11}&A^{jj'}_{12}\\A^{jj'}_{20}&A^{jj'}_{21}&A^{jj'}_{22}\end{array}\right].  $$
Explicitly, we obtain:
\begin{itemize}
\item  For $n=1,$ i.e. $x=x^1:$ 
$$\begin{array}{lll}
A^{jj'}_{00}=\frac{\partial^2L}{\partial u^j   \partial u^{j'}},
&A^{jj'}_{01}=\frac{\partial^2L}{\partial u^j   \partial u^{j'}_{x}},
&A^{jj'}_{02}=\frac{\partial^2L}{\partial u^j   \partial u^{j'}_{2x}},\\
A^{jj'}_{10}=\frac{\partial^2L}{\partial u^j_{x}   \partial u^{j'}},
&A^{jj'}_{11}=\frac{\partial^2L}{\partial u^j_{x}   \partial u^{j'}_{x}},
&A^{jj'}_{12}=\frac{\partial^2L}{\partial u^j_{x}   \partial u^{j'}_{2x}},\\
A^{jj'}_{20}=\frac{\partial^2L}{\partial u^j_{2x}   \partial u^{j'}},
&A^{jj'}_{21}=\frac{\partial^2L}{\partial u^j_{2x}   \partial u^{j'}_{x}},
&A^{jj'}_{22}=\frac{\partial^2L}{\partial u^j_{2x}   \partial u^{j'}_{2x}} , \end{array}$$
\item  For $n=2,$ i.e. $x=(x^1,x^2):$
$$\begin{array}{ll} A^{jj'}_{00}=\frac{\partial^2L}{\partial u^j   \partial u^{j'}},&
A^{jj'}_{01}=\left(\begin{array}{ll}\frac{\partial^2L}{\partial u^j   \partial u^{j'}_{x^1}}& \frac{\partial^2L}{\partial u^j   \partial u^{j'}_{x^2}}   \end{array}\right)\end{array},$$
$$\begin{array}{l}A^{jj'}_{02}=\left(\begin{array}{lll}\frac{\partial^2L}{\partial u^j   \partial u^{j'}_{2x^1}}&
\frac{\partial^2L}{\partial u^j   \partial u^{j'}_{x^1x^2}}&\frac{\partial^2L}{\partial u^j   \partial u^{j'}_{2x^2}}   \end{array}\right)\end{array}    ,$$
$$\begin{array}{ll}
A^{jj'}_{10}=\left(\begin{array}{l}\frac{\partial^2L}{\partial u^j_{x^1}   \partial u^{j'}}\\
 \frac{\partial^2L}{\partial u^j_{x^1}   \partial u^{j'}}\end{array}\right)  , &A^{jj'}_{11}=\left(\begin{array}{ll}\frac{\partial^2L}{\partial u^j_{x^1}   \partial u^{j'}_{x^1}}&
\frac{\partial^2L}{\partial u^j_{x^1}   \partial u^{j'}_{x^2}}\\
\frac{\partial^2L}{\partial u^j_{x^2}   \partial u^{j'}_{x^1}}&\frac{\partial^2L}{\partial u^j_{x^2}   \partial u^{j'}_{x^2}}         \end{array}\right)\end{array},$$
$$A^{jj'}_{12}=\left(\begin{array}{lll}
\frac{\partial^2L}{\partial u^j_{x^1}   \partial u^{j'}_{2x^1}}&\frac{\partial^2L}{\partial u^j_{x^1}   \partial u^{j'}_{x^1x^2}}&\frac{\partial^2L}{\partial u^j_{x^1}   \partial u^{j'}_{2x^2}}\\
\frac{\partial^2L}{\partial u^j_{x^2}   \partial u^{j'}_{2x^1}}&\frac{\partial^2L}{\partial u^j_{x^2}   \partial u^{j'}_{x^1x^2}}&\frac{\partial^2L}{\partial u^j_{x^2}   \partial u^{j'}_{2x^2}}
       \end{array}\right),$$
$$\begin{array}{ll}
A^{jj'}_{20}=\left(\begin{array}{l}\frac{\partial^2L}{\partial u^j_{2x^1}   \partial u^{j'}}\\
 \frac{\partial^2L}{\partial u^j_{x^1x^2}   \partial u^{j'}}\\
 \frac{\partial^2L}{\partial u^j_{2x^2}   \partial u^{j'}}\end{array}\right) , &
  A^{jj'}_{21}=\left(\begin{array}{ll}\frac{\partial^2L}{\partial u^j_{2x^1}   \partial u^{j'}_{x^1}}&
\frac{\partial^2L}{\partial u^j_{2x^1}   \partial u^{j'}_{x^2}}\\
\frac{\partial^2L}{\partial u^j_{x^1x^2}   \partial u^{j'}_{x^1}}&\frac{\partial^2L}{\partial u^j_{x^1x^2}   \partial u^{j'}_{x^2}}\\
\frac{\partial^2L}{\partial u^j_{2x^2}   \partial u^{j'}_{x^1}}&\frac{\partial^2L}{\partial u^j_{2x^2}   \partial u^{j'}_{x^2}}           \end{array}\right)\end{array},$$
$$A^{jj'}_{22}=\left(\begin{array}{lll}
\frac{\partial^2L}{\partial u^j_{2x^1}   \partial u^{j'}_{2x^1}}&\frac{\partial^2L}{\partial u^j_{2x^1}   \partial u^{j'}_{x^1x^2}}&\frac{\partial^2L}{\partial u^j_{2x^1}   \partial u^{j'}_{2x^2}}\\
\frac{\partial^2L}{\partial u^j_{x^1x^2}   \partial u^{j'}_{2x^1}}&\frac{\partial^2L}{\partial u^j_{x^1x^2}   \partial u^{j'}_{x^1x^2}}&\frac{\partial^2L}{\partial u^j_{x^1x^2}   \partial u^{j'}_{2x^2}}\\
\frac{\partial^2L}{\partial u^j_{2x^2}   \partial u^{j'}_{2x^1}}&\frac{\partial^2L}{\partial u^j_{2x^2}   \partial u^{j'}_{x^1x^2}}&\frac{\partial^2L}{\partial u^j_{2x^2}   \partial u^{j'}_{2x^2}}
       \end{array}\right).$$
\end{itemize}
\end{example}

Let $u=(u^1, \cdots  ,u^m)\in \left(\mathcal{C}^{s}(\Omega)\right)^m$ and $\phi=(\phi^1, \cdots  ,\phi^m)$ with $\phi^j\in \mathcal{C}^{\infty}(\Omega),$
the set of continuously infinitely differentiable functions on $\Omega.$
Consider the function $f$ defined by
$$ f(\epsilon)=J[u+\epsilon \phi]=\int_{\Omega}L\left(x,u^{(s)}(x)+\epsilon \phi^{(s)}(x)\right)dx,\quad \forall\,\epsilon >0.$$
Then, the first derivative of $f$ yields
\begin{eqnarray}
 f'(\epsilon)&=& \frac{d}{d \epsilon}J[u+\epsilon \phi]\nonumber\\
             &=& \int_{\Omega}\sum_{j=1}^{m}\sum_{k=0}^{s}\sum_{h=1}^{p_k}\,\phi^j_{(k)}[h](x)\,\frac{\partial L\left(x,u^{(s)}(x)
+\epsilon \phi^{(s)}(x)\right)}{\partial u^j_{(k)}[h]}\,dx^1\cdots dx^n\nonumber,
\end{eqnarray}
or equivalently
\begin{eqnarray}
  f'(\epsilon)&=&\int_{\Omega}\sum_{j=1}^{m}\,\sum_{k_1+\cdots+k_n=0}^{s}\,\,\phi^j_{k_1x^1\cdots k_nx^n}(x)\,\frac{\partial L\left(x,u^{(s)}(x)+\epsilon \phi^{(s)}(x)\right)}{\partial u^j_{k_1x^1\cdots k_n x^n}}\,dx^1\cdots dx^n.\nonumber
 \end{eqnarray}
The second derivative of $f$ is given by
\begin{eqnarray}
f''(\epsilon)&=& \frac{d^2}{d\epsilon^2}J[u+\epsilon \phi]\nonumber\\
             &=&\int_{\Omega}\sum_{j,j'=1}^m \sum_{k,k'=0}^s \sum_{h=1}^{p_k}\sum_{h'=1}^{p_{k'}}\phi^j_{(k)}[h](x)\phi^{j'}_{(k')}[h'](x) \frac{\partial^2 L\left(x,u^{(s)}(x)+\epsilon \phi^{(s)}(x)\right)}{\partial u^j_{(k)}[h]\,\partial u^{j'}_{(k')}[h']}dx\nonumber\\
            & = & \int_{\Omega}\sum_{j,j'=1}^m \sum_{k,k'=0}^s \,\phi^j_{(k)}(x)\,  A^{jj'}_{kk'}\left(x,u^{(s)}(x)+\epsilon \phi^{(s)}(x)\right)  \,^{t}\phi^{j'}_{(k')}(x)\,dx\nonumber\\
            & = & \int_{\Omega}\sum_{j,j'=1}^m\,\phi^{j(s)}(x)\, A^{jj'}\left(x,u^{(s)}(x)+\epsilon \phi^{(s)}(x)\right)  \,^{t}\phi^{j'(s)}\,dx\nonumber\\
            & = & \int_{\Omega}\phi^{(s)}(x)\, A\left(x,u^{(s)}(x)+\epsilon \phi^{(s)}(x)\right) \, ^{t}\phi^{(s)}(x)\,dx\nonumber
\end{eqnarray}
leading to the simplified form
\begin{eqnarray}
  f''(\epsilon)  &=&\int_{\Omega}\sum_{j,j'=1}^{m}\,\sum_{k_1+\cdots +k_n=0}^{s}\,\, \sum_{k'_1+\cdots +k'_n=0}^{s}\frac{\partial^2 L\left(x,u^{(s)}(x)+\epsilon \phi^{(s)}(x)\right)}{\partial u^j_{k_1x^1\cdots k_nx^n}\partial u^{j'}_{k'_1x^1\cdots k'_nx^n}}^{}\nonumber\\
   & \times&
  \phi^j_{k_1x^1\cdots k_nx^n}(x)\, \phi^{j'}_{k'_1x^1\cdots k'_nx^n}(x) \,   dx.\nonumber
 \end{eqnarray}
Provided  all these definitions, we  may now state the  main results of this work.
\begin{theorem}\emph{\textbf{(Sufficient condition of extrema)}}\\
 Let $W$ be an open subspace of the Fr\'{e}chet space $\left(\mathcal{C}^{s}(\Omega)\right)^m$ and $J:W\rightarrow \mathbb{R}$ a
 functional defined by (\ref{eqn0}), continuously twice differentiable at $\overline{u}\in W$ such that
\begin{equation}
 \frac{\delta}{\delta u^j} L\left(x,\overline{u}^{(s)}(x)\right)  =0,\quad j=1,\cdots,m,
\end{equation}
i.e. the function $\overline{u}$ is a critical point of the functional  $J.$
Then, $J$ admits a local minimum (resp. maximum)  at the critical point $\overline{u}$ if
the corresponding matrix $A\left(x,\overline{u}^{(s)}(x)\right)$ defined by (\ref{eqac}) is positive (resp. negative) definite for all $x\in \Omega.$
 \end{theorem}
 \begin{prof}
 For all $\phi=(\phi^1, \cdots  ,\phi^m)$ with $\phi^j\in \mathcal{C}^{\infty}(\Omega),$  the
 second derivative of the  real valued function
 $f$, defined by $ f(\epsilon)=J[\overline{u}+\epsilon \phi],$ at $\epsilon=0$ is
 \begin{eqnarray}
  f''(0)&=&\frac{d^2}{d\epsilon^2}J[\overline{u}+\epsilon \phi]|_{\epsilon =0}\nonumber\\
&=&  \int_{\Omega}\phi^{(s)}(x)\, A\left(x,\overline{u}^{(s)}(x)\right) \, ^{t}\phi^{(s)}(x)\,dx.\nonumber
 \end{eqnarray}
 Thus,  if the  matrix $A\left(x,\overline{u}^{(s)}(x)\right)$ defined by (\ref{eqac}) is positive (resp. negative) definite for all $x\in \Omega,$ 
then $J''(\overline{u})=f''(0)>0\,\,(\emph{\emph{resp.}} < 0),$
 i.e. the function $\overline{u}$ is a local minimum (resp. maximum) point for the functional $J.$
 \end{prof}

 \texttt{ }

 Assume now that $\Omega = \prod_{i=1}^{n}]a_i,\,b_i  [$ with $a_i,\,b_i\in \mathbb{R}.$
 \begin{theorem}\emph{\textbf{(Necessary condition of extrema)}}\\
 Let $W$ be an open subspace of the Fr\'{e}chet space $\left(\mathcal{C}^{s}(\Omega)\right)^m$ and $J:W\rightarrow \mathbb{R}$ a
 functional defined by (\ref{eqn0}),  continuously twice differentiable on $W.$
Let $\overline{u}\in W$ be  a local minimum (resp. maximum) point of $J.$
Then for all $x\in \Omega,$
\begin{equation}
 B_l(x)\equiv\sum_{j,j'=1}^{m}\sum_{h,h'=1}^{p_l}\frac{\partial^2 L\left(x,\overline{u}^{(s)}(x)\right)}{\partial u^j_{(l)}[h]  \partial u^{j'}_{(l)}[h']}\,   \geq 0 \,\,(\emph{\emph{resp.}}   \leq  0),\quad l=0,1,2,\cdots,s.
\end{equation}
 \end{theorem}
 \begin{prof}
Let $\overline{x}=(\overline{x}^1,\cdots  ,\overline{x}^n)\in \Omega.$ Then there exists $r_0>0$ such that the ball of radius $r_0$ centered at $\overline{x},$
$B(\overline{x};{r_0})\subset \Omega.$ Choosing $0<\epsilon<\frac{r_0}{2n}$ implies
$\prod_{i=1}^n\left[\overline{x}^i-\epsilon,\,\overline{x}^i+\epsilon\right]\subset B(\overline{x};{r_0}).$
For $l=0,1,2,\cdots  ,s$ we define some particular functions $\phi_l$  by
\begin{equation}\label{eq0c1}
\phi_l(x)=\sum_{i=1}^n \phi_{l,i}(x^i),\quad \forall\,x\in \prod_{i=1}^n\left[\overline{x}^i-\epsilon,\,\overline{x}^i+\epsilon\right],
\end{equation}
where for $i=1,2,  \cdots  ,n$
\begin{equation}
 \phi_{l,i}(x^i)=\left\{\begin{array}{cl}
 0, & a_i\leq x^i\leq \overline{x}^i-\epsilon\\
 1-\frac{1}{\epsilon}(x^i-\overline{x}^i)^l\,\emph{\emph{sign}}(x^i-\overline{x}^i),& \overline{x}^i-\epsilon\leq x^i\leq \overline{x}^i+\epsilon\\
 0,& \overline{x}^i+\epsilon\leq x^i\leq b_i.
        \end{array}\right.
\end{equation}
Clearly, the functions $\phi_l$ are continuously infinitely differentiable everywhere except at the end points of its domain.
Consider the real valued function $f$
defined by   
\begin{eqnarray}
 f(\widetilde{\epsilon})&=& \widetilde{J}[\overline{u}+\widetilde{\epsilon} \phi]\nonumber\\
 &=& \int_{\prod_{i=1}^n\left[\overline{x}^i-\epsilon,\,\overline{x}^i+\epsilon\right]}L\left(x,\overline{u}^{(s)}(x)+\widetilde{\epsilon} \phi^{(s)}(x)\right)dx,\quad \forall \, \widetilde{\epsilon} > 0,
\end{eqnarray}
with $\phi=(\phi^1, \cdots  ,\phi^m),$ where the real valued functions
  $\phi^j$ are continuously infinitely differentiable everywhere except 
at the end points of the domain  $\prod_{i=1}^n\left[\overline{x}^i-\epsilon,\,\overline{x}^i+\epsilon\right].$  Then,
\begin{eqnarray}\label{eq0c2}
  f''(0)&=&\frac{d^2}{d\widetilde{\epsilon}^2}J[\overline{u}+\widetilde{\epsilon} \phi]|_{\widetilde{\epsilon} =0}\nonumber\\
  &=& \int_{\prod_{i=1}^n\left[\overline{x}^i-\epsilon,\,\overline{x}^i+\epsilon\right]}\sum_{j,j'=1}^m \sum_{k,k'=0}^s \sum_{h=1}^{p_k}\sum_{h'=1}^{p_{k'}} \frac{\partial^2 L\left(x,\overline{u}^{(s)}(x)\right)}{\partial u^j_{(k)}[h]\,\partial u^{j'}_{(k')}[h']}\nonumber\\
  & \times & \phi^j_{(k)}[h](x)\phi^{j'}_{(k')}[h'](x) dx. \label{eq0c2}
 \end{eqnarray}
Substitute in (\ref{eq0c2})  the functions $\phi^j$ and $\phi^{j'},\,\, j,j'=1,2, \cdots ,m$ by the same function  $\phi_l$ defined by (\ref{eq0c1}) and let $\epsilon\rightarrow 0.$ We obtain, after computation and simplification,
\begin{eqnarray}
f''(0)&=&\lim_{\epsilon\rightarrow 0}\left( \int_{\prod_{i=1}^n\left[\overline{x}^i-\epsilon,\,\overline{x}^i+\epsilon\right]}\sum_{j,j'=1}^m  \sum_{h=1}^{p_l}\sum_{h'=1}^{p_{l}}\phi^j_{(l)}[h](x)\phi^{j'}_{(l)}[h'](x) \frac{\partial^2 L\left(x,\overline{u}^{(s)}(x)\right)}{\partial u^j_{(k)}[h]\,\partial u^{j'}_{(k')}[h']}dx\right)\nonumber\\
&=&\lim_{\epsilon\rightarrow 0}\left( \frac{1}{\epsilon^2}\int_{\prod_{i=1}^n\left[\overline{x}^i-\epsilon,\,\overline{x}^i+\epsilon\right]}B_l(x)\,\,dx^1\cdots dx^n\right). \label{eqcs}
\end{eqnarray}
Since the function $\overline{u}$ is a local minimum (resp. maximum) for the functional $J,$ then
we must have
\begin{equation}\label{eqcs2}
f''(0)=J''(\overline{u}) \geq 0 \quad (\emph{\emph{resp.}} \leq 0).
\end{equation}
Taking into account the equality (\ref{eqcs}), the  conditions (\ref{eqcs2}) lead to the inequalities
\begin{equation}
 B_l(\overline{x}) \geq 0 \quad (\emph{\emph{resp.}} \,\,  \leq  0).
 \end{equation}
This end the proof since $\overline{x}$ is an arbitrary point in the open subset $\Omega.$
 \end{prof}
 \\



\section*{Acknowledgments}
This work is partially supported by the ICTP through the
OEA-ICMPA-Prj-15. The ICMPA is in partnership with
the Daniel Iagolnitzer Foundation (DIF), France.



\end{document}